\documentclass[a4paper,12pt]{article}

\usepackage[frenchb,english]{babel}
\usepackage{amsmath,amssymb,amscd,amsfonts}
\usepackage[T1]{fontenc}
\usepackage[latin1]{inputenc}
\usepackage{indentfirst}
\usepackage[all,ps]{xy}

\newcommand{\cs}{contact sphere}
\newcommand{\cc}{contact circle}
\newcommand{\cps}{contact $p$-sphere}
\newcommand{\css}{contact sphere }
\newcommand{\cpss}{contact $p$-sphere }
\newcommand{\ccs}{contact circle }
\newcommand{\sph}{\mathbb{S}}
\newcommand{\tor}{\mathbb{T}}
\newcommand{\ha}{\mathbb{H}}
\newcommand{\R}{\mathbb{R}}

\newcommand{\Sc}{\mathcal{S}_c}

\newcommand{\sla}{\sum_{i=1}^3 \lambda_i^2}
\newcommand{\slas}{\sum_{i=1}^3 \lambda_i}
\newcommand{\finthm}{\nopagebreak{} \begin{flushright}
$\blacksquare$
\end{flushright}}
\newcommand{\finlem}{\nopagebreak{} \begin{flushright}
$\square$
\end{flushright}}
\newcommand{\dop}{~: }

\newtheorem{defi}{Definition}
\newtheorem{theorem}{Theorem}
\newtheorem{cor}{Corollary}
\newtheorem{lemma}{Lemma}
\newtheorem{prop}{Proposition}

\begin{document} 
\sloppy
\author{Mathias Zessin}

\title{On contact $p$-spheres}

\date{}
\maketitle
\vspace{-1cm}

{\small \begin{center}
        Laboratoire de Mathématiques \\
        Université de Mulhouse \\
        4, rue des frères Lumière \\
        68093 Mulhouse Cedex \\
        \underline{e-mail:} Mathias.Zessin@uha.fr
        \end{center}}
\begin{abstract}
We study invariant \cps s on principal $\sph^1$-bundles and solve the corresponding existence problem in dimension 3. 
Moreover, we show that contact $p$-spheres can only exist on $(4n-1)$-dimensional manifolds and we construct examples of \cps s on such manifolds. We also consider relations between tautness and roundness, a regularity property concerning the Reeb vector fields of the contact forms in a \cps.
\end{abstract}

\section{Introduction and generalities} 

Contact circles and contact $p$-spheres are families of contact forms parametrized by the circle and the $p$-sphere respectively and have been introduced by H. Geiges and J. Gonzalo in 1995 (see \cite{GG1}). These authors give a fundamental existence theorem for \cc s and closely study a special class of \cc s with additional geometrical properties (see also \cite{GG2}). Their work is about 3-dimensional manifolds.

In the first part of this paper we are interested in \emph{invariant} \cc s and contact $p$-spheres on principal $\sph^1$-bundles. In this specific situation we adapt some methods developed by R. Lutz which have been used to construct the first examples of contact forms on the 5-dimensional torus (see \cite{L2}). They will be particularly helpful for the construction of new classes of examples. In this work, we will see examples of invariant \cc s and \cs s on 3-dimensional principal $\sph^1$-bundles. In the second part, we examine the general situation in higher dimensions and give examples of contact $p$-spheres there.

In paragraph 2 we start with some examples. In dimension 3, the most natural ones are found on $\sph^3$. In paragraph 3, we develop tools to study invariant \cc s and \cs s. In paragraphs 4 and 5 we state the main theorems of this first part. They solve the existence problem of invariant \cc s and \cs s on principal $\sph^1$-bundles over closed, connected and orientable surfaces. In paragraph 6 we construct examples on all such principal $\sph^1$-bundles where \cc s might exist according to the theorems of paragraphs 4 and 5. In paragraph 7 we consider higher dimensions and prove that on manifolds of dimension $4n+1$, \cc s and contact $p$-spheres do not exist. In paragraph 8 we study round contact $p$-spheres, that is, contact $p$-spheres whose Reeb vector fields are stable under linear combinations within the contact $p$-sphere, and in paragraph 9 we give examples of contact $p$-spheres on the spheres of dimension $4n-1$. In paragraph 10 we make a connection to Sasakian geometry, which gives more examples of contact spheres.

\vspace{0.4cm}
Let us now give some definitions and preliminary results.

Let $M$ be a $(2n+1)$-dimensional manifold. A contact form on $M$ is a 1-form $\omega$ such that $\omega \wedge (d\omega)^n$ is a volume form on $M$. The Reeb vector field $R$ associated to $\omega$ is defined by $\omega(R)=1$ and $R\lrcorner \, d\omega=0$.

\begin{defi} ---
The set of normalized linear combinations of p+1 contact forms $\omega_1,\dots,\omega_{p+1}$ 
$$ \Big\{ \quad \sum_{i=1}^{p+1} \lambda_i \,\omega_i \quad | \quad \sum_{i=1}^{p+1} \lambda_i^2=1 \quad \Big\} $$
is called a \emph{contact $p$-sphere}, if every element of this family is a contact form. It will be denoted by $\mathcal{S}^p_c\{\omega_1,\dots,\omega_{p+1}\}$.
\end{defi}

Note that Geiges and Gonzalo call \ccs a pair $(\omega_1,\,\omega_2)$ of contact forms whose normalized linear combinations are still contact forms and \css a triple $(\omega_1,\,\omega_2,\,\omega_3)$ of such forms. Here, however, contact $p$-sphere will denote the whole family generated by these forms.

In this paper, a contact 1-sphere will be called a \emph{\cc}, and a contact 2-sphere might simply be called \emph{\cs}. Note that each \css contains \cc s and that in general each contact $m$-sphere contains contact $n$-spheres for $n \leq m$.

Following Geiges and Gonzalo, a \ccs or a contact $p$-sphere will be called \emph{taut}, if all its elements generate the same volume form, that is, if for any element $\omega=\lambda_1\,\omega_1+\lambda_2\,\omega_2$ of, say, a \ccs in dimension 3, the form $\omega \wedge d\omega$ does not depend on its coefficients $\lambda_1$ and $\lambda_2$.

\section{Examples on $\sph^3$ and on $\tor^3$} 
\label{exemples}
We use a quaternionic representation of $\sph^3$ to get the following natural example of a \cs\dop

In the space of quaternions $\ha$, the forms
\begin{eqnarray*}
\left \{ \begin{array}{rcl}
\tilde{\alpha}&=&<qi,\,dq>\quad=\quad q_1\,dq_2-q_2\,dq_1+q_4\,dq_3-q_3\,dq_4\\[1mm]
\tilde{\beta}&=&<qj,\,dq> \quad = \quad q_1\, dq_3-q_3\,dq_1+q_2\,dq_4-q_4\,dq_2\\[1mm]
\tilde{\gamma}&=&<qk,\,dq> \quad = \quad q_3\, dq_2-q_2\,dq_3+q_1\,dq_4-q_4\,dq_1
\end{array} \right.
\end{eqnarray*}
induce three contact forms $\alpha,\beta$ and $\gamma$ on the unit sphere $\sph^3$, which are linearly independent and complementary in the sense that their Reeb vector fields give an orthonormal basis of the tangent space $T_p \sph^3$ at every point $p \in \sph^3$ with respect to the induced Euclidean metric.

Any form 
$$
\tilde{\omega}=\lambda_1\,\tilde{\alpha}+\lambda_2\,\tilde{\beta}+\lambda_3\,\tilde{\gamma} \quad\textrm{with}\quad\sla=1
$$
satisfies
$$
\tilde{\omega}\wedge\,d\tilde{\omega}\wedge\,(q_1\,dq_1+q_2\,dq_2+q_3\,dq_3+q_4\,dq_4)=dq_1\wedge \,dq_2\wedge \,dq_3\wedge \,dq_4,
$$
which is non-zero everywhere. Hence the induced form $\omega=\lambda_1\,\alpha+\lambda_2\,\beta+\lambda_3\,\gamma$ is a contact form and $\omega \wedge d\omega=\alpha \wedge d\alpha$. This gives a taut \cs. It is also round in a sense defined in paragraph \ref{round cs}.

We also observe that the forms $\alpha,\,\beta$ and $\gamma$ are invariant with respect to the vector field $X$ induced by  $$q_1\,\frac{\partial}{\partial{q_2}}-q_2\,\frac{\partial}{\partial{q_1}}+q_3\,\frac{\partial}{\partial{q_4}}-q_4\,\frac{\partial}{\partial{q_3}},$$
which is the Reeb vector field of the contact form $\hat{\alpha}$ induced by $<iq,\,dq>$ on $\sph^3$. The orbits of $X$ are the fibres of a principal Hopf fibration of $\sph^3$ over $\sph^2$ with fibre $\sph^1$ and connection form $\hat{\alpha}$. In this case, $\Sc^2\{\alpha,\beta,\gamma\}$ is an invariant \css on this $\sph^1$-bundle.

\begin{displaymath}
\xymatrix{
\quad \quad \; \; \sph^3 \ar[d]_{\pi(q)=q\,i\,\bar{q}}^{\; \sph^1} \subset \ha \\
\sph^2
}
\end{displaymath}

\vspace{0.4cm}
Another quite natural example of \cc s can be found on $\tor^3$, with pseudo-coordinates $(\theta_1,\,\theta_2,\,\theta_3)$. On this torus, consider the forms

\begin{eqnarray*}
\left \{ \begin{array}{rcl}
\omega_1 &=& \phantom{-} \cos(n\,\theta_1) d\theta_2+\sin(n\,\theta_1) d\theta_3 \\
\omega_2 &=& -\sin(n\,\theta_1) d\theta_2+\cos(n\,\theta_1) d\theta_3,
\end{array} \right.
\end{eqnarray*}
for some non-zero integer $n$. Setting $\omega=\lambda_1\,\omega_1+\lambda_2\,\omega_2$, with $\lambda_1^2+\lambda_2^2=1$, we get
$$\omega \wedge d\omega=-n \, d\theta_1 \wedge d\theta_2 \wedge d\theta_3.$$ 
Thus for any $n$, $\omega_1$ and $ \omega_2$ generate a taut \cc. We can observe that $\omega_1$ and $\omega_2$ and thus all elements of the \ccs are invariant with respect to the vector fields $\frac{\partial}{\partial \theta_2}$ and $\frac{\partial}{\partial \theta_3}$. So these \cc s are even invariant with respect to the corresponding $\tor^2$-action.

\section{Singular sets and knotted fibrations} \label{ens sing}

In this paragraph we develop tools to study invariant \cc s and \cps s on principal $\sph^1$-bundles over connected closed manifolds. We first examine singular sets associated to invariant contact forms in this context and then consider these sets as knots of knotted fibrations. This approach allows great insight into the topology of invariant \cc s and \cps s.

\vspace{0.4cm} 
Consider a $(2n+1)$-dimensional manifold $M$, which is a principal $\sph^1$-bundle over a $(2n)$-dimensional manifold $B$ with connection form $\alpha$. Assume that $B$ is a compact, connected, orientable manifold without boundary. Then for each invariant form $\omega$ on $M$, there is a form $\eta$ and a function $\varphi$ on the base space $B$, such that 
$$ \omega=\pi^*(\eta)+\pi^*(\varphi)\,\alpha, $$
where $\pi$ is the bundle map.

We have the following technical lemma\dop

\begin{lemma} \label{lem fond} --- Let $\mathcal{S}^p_c\{\omega_1,\dots,\omega_{p+1}\}$ be an invariant contact $p$-sphere on $M$ and assume that for each $i$, $$\omega_i=\pi^*(\eta_i)+\pi^*(\varphi_i)\,\alpha $$
for some forms and functions $\eta_i$ and $\varphi_i$ on $B$. Then for any coefficients  $(\lambda_1,\dots,\lambda_{p+1})$ with $\sum_{i=1}^{p+1} \lambda_i^2=1$, the expressions
\begin{equation} \label{func-0}
\sum_{i=1}^{p+1} \lambda_i\,\varphi_i \textrm{ and } \sum_{i=1}^{p+1} \lambda_i\,d\varphi_i
\end{equation}
do not vanish simultaneously at any point of $B$.
\end{lemma}

{\bf Proof\dop} --- This is an immediate consequence of the contact property of linear combinations of the forms $\omega_1,\dots \omega_{p+1}$. Indeed, for $\omega=\sum_{i=1}^{p+1} \lambda_i\,\omega_i$, we have 
$$ \omega\wedge (d\omega)^n=(\sum_{i=1}^{p+1}\lambda_i\,\eta_i +\sum_{i=1}^{p+1}\lambda_i\,\varphi_i\,\alpha)\wedge (\sum_{i=1}^{p+1}\lambda_i\,d\eta_i+\sum_{i=1}^{p+1}\lambda_i\,d\varphi_i\wedge\alpha+\sum_{i=1}^{p+1}\lambda_i\,\varphi_i\,d\alpha)^n.$$
This term equals 
$$ (\sum_{i=1}^{p+1}\lambda_i\,\eta_i)\wedge (\sum_{i=1}^{p+1}\lambda_i\,d\eta_i)^n,$$
at points where the expressions (\ref{func-0}) vanish simultaneously. But a $(2n+1)$-form on $B_{2n}$ is zero. This contradicts the assumption that $\mathcal{S}^p_c\{\omega_1,\dots,\omega_{p+1}\}$ is a contact $p$-sphere and proves the lemma.
\finlem

Let us now define on the base manifold the singular set associated to a given invariant contact form on $M$, following R. Lutz (see \cite{L1}).

\begin{defi} ---
Let $\omega=\pi^*(\eta)+\pi^*(\varphi)\,\alpha$ be an invariant contact form on an $\sph^1$-bundle $M$ over a manifold $B$. The \emph{singular set} associated to $\omega$ is the set $\Sigma_\omega =\varphi^{-1}(0)$ on $B$.
\end{defi}

We now prove some fundamental properties of the singular sets of invariant \cc s and \cps s on a principal $\sph^1$-bundle $M$ over a closed connected manifold $B$.

\begin{prop} ---
\begin{description}
\item[{\it i)}] The singular sets are submanifolds of $B$.
\item[{\it ii)}] If $M$ admits an invariant \cc, then each point of $B$ lies on the singular set of some element of this \cc.
\item[{\it iii)}] The singular sets of two forms of the same \ccs or \cpss are isotopic.
\item[{\it iv)}] The singular set of an element of an invariant \ccs or \cpss is non-empty.
\item[{\it v)}] The singular sets of two different and non-opposite elements of an invariant \ccs or \cpss intersect transversally.
\end{description}
\end{prop}

{\bf Proof\dop} ---
According to R. Lutz (\cite{L1}) (or to Lemma \ref{lem fond} in the particular case of $p=0$), $d\varphi$ does not vanish on $\varphi^{-1}(0)$, if $\omega=\pi^*(\eta)+\pi^*(\varphi)\,\alpha$ is an invariant contact form on $M$. Thus the corresponding singular set is a submanifold of $B$.

Let $\omega_1$ and $\omega_2$ generate an invariant \ccs on $M$. Writing, for $i=1,2$,
$$ \omega_i=\pi^*(\eta_i)+\pi^*(\varphi_i)\,\alpha,$$
an arbitrary element of the \ccs will be 
\begin{eqnarray*}
\omega&=&\sin(\theta)\,\omega_1+\cos(\theta)\,\omega_2\\ &=&\pi^*(\sin(\theta)\,\eta_1+\cos(\theta)\,\eta_2)+\pi^*((\sin(\theta)\,\varphi_1+\cos(\theta)\,\varphi_2))\,\alpha\\
&=&\pi^*(\eta)+\pi^*(\varphi)\,\alpha.
\end{eqnarray*}
Thus, for any values $\varphi_1(p)$ and $\varphi_2(p)$, there is some $\theta\in [0,2\pi]$ such that $\varphi(p)=0$.

This argument shows also that in any \ccs there are elements whose singular sets are not empty. On the other hand, according to \cite{Gr} (see also \cite{Ma}), the contact structures associated to elements of one-parameter-families of contact forms are isomorphic. Obviously, \cc s are special one-parameter-families, so the corresponding structures are isomorphic. The isomorphisms are equivariant, as explained in \cite{L1}, so the singular sets of elements of a \ccs are isotopic. Thus any element of a \ccs has a non-empty singular set.

 In contact $p$-spheres, $p$ being an arbitrary dimension, any two different and non-opposite elements generate a \cc, so their singular sets are non-empty and isotopic, too.
 
 Finally, at a point where the singular sets of two different and non-opposite elements $\omega_1=\pi^*(\eta_1)+\pi^*(\varphi_1)\,\alpha$ and $\omega_2=\pi^*(\eta_2)+\pi^*(\varphi_2)\,\alpha$ of an invariant \ccs intersect, we have $\varphi_1=\varphi_2=0$, so by Lemma \ref{lem fond}, the forms $d\varphi_1$ and $d\varphi_2$ are linearly independent, hence the intersection of $\Sigma_{\omega_1}$ and $\Sigma_{\omega_2}$ is transversal.
\finlem

\vspace{0.4cm}
One way to study singular sets and thus invariant contact $p$-spheres is to consider their associated knotted fibrations. This viewpoint is due to R. Lutz and generalizes J. Milnor's open book constructions in the context of invariant contact forms. It means the following (see \cite{L2})\dop

\begin{defi} ---
$(B,\varphi,N)$ is a knotted fibration along $N$ over $\sph^{k-1}$, if 
\begin{enumerate}
	\item $B$ is a connected, orientable, compact manifold
	\item $N$ is either empty or a closed codimension-k-submanifold in $B$
	\item $\varphi:B\setminus N\longrightarrow \sph^{k-1}$ is a locally trivial fibration
	\item If $N$ is non-empty, there exists an open neighborhood $W$ of $N$ and a diffeomorphism $h:N\times D^k\longrightarrow W$, where $D^k$ is the unit disk of $\R^k$, such that $h(z,0)=z$ on $N\times \{0\}$ and such that the following diagram commutes\dop
\begin{displaymath}
\xymatrix{
N\times \left( D^k\setminus\{0\}\right) \ar[r]^{\quad \quad h} \ar[d]_{proj.}
&W\setminus N \ar[d]_{\varphi}
\\
D^k\setminus\{0\} \ar[r]^{\frac{\centerdot}{\|\centerdot\|}}
&\sph^{k-1}
}
\end{displaymath}
\end{enumerate}
\end{defi}

$N$ will be called the knot and $\varphi^{-1}(p)$ will be called a fibre of the knotted fibration, for any point $p \in \sph^{k-1}$.

\vspace{0.4cm}
As an example in dimension 3, we can consider singular sets associated to invariant \cps s on principal $\sph^1$-bundles over surfaces in different ways\dop
\begin{itemize}
	\item Given one contact form $\omega=\pi^*(\eta)+\pi^*(\varphi)\,\alpha$, the singular set $\Sigma_\omega =\varphi^{-1}(0)$ is a curve on the base space and defines a knotted fibration over $\sph^0$ with knot $\Sigma_\omega$ and two fibres determined by the sign of the function $\varphi$. These fibres are well defined (see \cite{L1}).
	
	\item For a \ccs generated by two invariant contact forms $\omega_1$ and $\omega_2$, the intersection $\Sigma$ of the associated singular sets $\Sigma_{\omega_1}$ and $\Sigma_{\omega_2}$ is finite. $\Sigma$ is the knot of a knotted fibration of $B$ over $\sph^1$, whose fibres are curves starting and ending at points of $\Sigma$.

	\item In the case of an invariant \css with generating forms $\omega_1,\,\omega_2$ and $\omega_3$, there is no common singular set, which is due to Lemma \ref{lem fond}, but we can consider an associated knotted fibration of $M$ over $\sph^2$ with empty knot and fibres which are finite subsets of $B$. 
\end{itemize}

\vspace{0.3cm}
These are special cases of the following fibration theorem which gives this decomposition in full generality\dop

\begin{theorem} --- Fibration theorem \label{thm fibr}

Let $M$ be a (2n+1)-dimensional principal $\sph^1$-bundle over a closed connected base manifold $B_{2n}$ with connection form $\alpha$ and let $\mathcal{S}^p_c\{\omega_1,\dots,\omega_{p+1}\}$ be an invariant contact $p$-sphere on $M$. Write $\omega_i$ as $\omega_i=\pi^*(\eta_i)+\pi^*(\varphi_i)\,\alpha$, for i=1,\dots,p+1, and let  $\Sigma=\bigcap_{i=1}^{p+1} \varphi_i^{-1}(0)$ be the intersection of the singular sets of the generating forms. Then for $r=\sqrt{\sum_{i=1}^{p+1}\varphi_i^2}$, $$(B,\,(\frac{\varphi_1}{r},\dots,\frac{\varphi_{p+1}}{r}),\Sigma)$$ 
defines a knotted fibration over $\sph^p$.
\end{theorem}

{\bf Proof\dop} According to \cite{L2}, it is enough to show that the map $$\Phi=(\frac{\varphi_1}{r},\dots,\frac{\varphi_{p+1}}{r}):B\setminus\Sigma \rightarrow \sph^p$$
is of rank $p$ everywhere and that at any point of $\Sigma$, we have
$$ d\varphi_1\wedge\dots\wedge d\varphi_{p+1}\neq0.$$

The first condition is a consequence of Lemma \ref{lem fond}. Indeed, we show that if $\Phi$ is not of rank $p$, then there are coefficients $(\lambda_1,\,\dots,\,\lambda_{p+1})$ with $\sum_{i=1}^{p+1} \lambda_i^2=1$, such that
$$ \sum_{i=1}^{p+1}\lambda_i\varphi_i=0 \textrm{ and } \sum_{i=1}^{p+1}\lambda_i\,d\varphi_i=0, $$
at some points of $B\setminus\Sigma$, which is excluded by Lemma \ref{lem fond}.

On a point $x \in B\setminus\Sigma$, at least one of the functions $\varphi_i$ does not vanish. Assume that $\varphi_{p+1}$ is non-zero. Then the rank of $\Phi$ is given by the rank of the system $(d(\frac{\varphi_1}{r}),\dots,d(\frac{\varphi_p}{r}))$.

Assume now that there is some point $x\in B\setminus\Sigma$ and $p$ real numbers $\mu_1,\dots,\mu_p$, such that $\sum_{i=1}^p\mu_id(\frac{\varphi_i}{r})=0$. Differentiating this and using $dr=\frac{1}{r}(\sum_{i=1}^{p+1}\varphi_i d\varphi_i)$, we get
$$0=\frac{1}{r}(\sum_{i=1}^p \mu_id\varphi_i)-\sum_{i=1}^p\frac{\mu_i\varphi_i}{r^2}dr=\sum_{i=1}^{p+1} \lambda_i\, d\varphi_i, $$
where we set $\lambda_i=\frac{\mu_i}{r}-(\sum_{j=1}^p\frac{\mu_j\varphi_j}{r^2})\frac{\varphi_i}{r}$, for $i=1,\dots,p$ and $\lambda_{p+1}=-(\sum_{j=1}^p\frac{\mu_j\varphi_j}{r^2})\frac{1}{r}\varphi_{p+1}$. With the same coefficients, we have now
\begin{eqnarray*}
\sum_{i=1}^{p+1}\lambda_i\varphi_i&=&\sum_{i=1}^p\frac{\mu_i\varphi_i}{r}-\sum_{i=1}^p\varphi_i^2\sum_{j=1}^p\frac{\mu_j\,\varphi_j}{r^3}-\sum_{i=1}^p\frac{\mu_i\,\varphi_i}{r^3}\varphi_{p+1}^2\\
&=&\sum_{i=1}^p\frac{\mu_i\,\varphi_i}{r^3}(r^2-\sum_{j=1}^{p+1}\varphi_j^2)\\
&=&0.
\end{eqnarray*}

The second condition is another consequence of Lemma \ref{lem fond}. As on $\Sigma$ the functions $\varphi_i$ vanish, any expression $\sum_{i=1}^{p+1} \lambda_i\,d\varphi_i$ is non-zero by Lemma \ref{lem fond}. Thus the forms $\{d\varphi_i,\ i=1\dots p+1\}$ are linearly independent on $\Sigma$, hence $d\varphi_1 \wedge \dots \wedge d\varphi_{p+1} \neq 0$ on $\Sigma$. This completes the proof.

\finthm

\begin{cor} \label{S1 fibr} ---
Let $\Sc^1\{\omega_1,\,\omega_2\}$ be an invariant \ccs on an $\sph^1$-bundle $M$ over a closed connected surface $B$ with common singular set $\Sigma$. Then $B \setminus \Sigma$ fibres over $\sph^1$.
\end{cor}

This corollary is a useful tool in certain situations. It gives us some additional information about the singular sets associated to invariant \cc s on principal $\sph^1$-bundles over surfaces, for example the following properties.

\begin{prop} --- \label{descr cc}
\begin{description}
\item[a)] The singular set $\Sigma_\omega$ of an element $\omega$ of a \ccs generated by $\omega_1$ and $\omega_2$ on a principal $\sph^1$-bundle over $\sph^2$ is a topological circle and the singular sets of two different and non-opposite elements intersect in two points.
\item[b)] The singular set $\Sigma_\omega$ of an element $\omega$ of a \ccs generated by $\omega_1$ and $\omega_2$ on a principal $\sph^1$-bundle over $\tor^2$ is the union of an even number of topological circles and the singular sets of two different and non-opposite elements do not intersect.
\end{description}
\end{prop}

{\bf Proof\dop} In the case $B=\sph^2$, it is enough to show that $\Sigma_{\omega_1}$ is a topological circle. As it is a non-empty submanifold of $\sph^2$, it is a union of topological circles. To see why $\Sigma_{\omega_1}$ can not have several connected components, we first observe that the common singular set $\Sigma=\Sigma_{\omega_1} \cap \Sigma_{\omega_2}$ is the union of two points, as $\sph^2 \setminus \Sigma$ fibres over $\sph^1$, and that the fibres are non-closed curves starting and ending at different components of $\Sigma$. As the fibration map is $\Phi=(\frac{\varphi_1}{r},\,\frac{\varphi_2}{r})$, where $r=\sqrt{\varphi_1^2+\varphi_2^2}$, $\, \Sigma_{\omega_1} \setminus \Sigma$ is the union of two fibres which are given by $\Phi^{-1}(0,1)$ and $\Phi^{-1}(0,-1)$ and which meet on $\Sigma$. Thus $\Sigma_{\omega_1}$ is a topological circle.  

In the case $B=\tor^2$, the common singular set $\Sigma$ is empty, as $\tor^2 \setminus \Sigma$ fibres over $\sph^1$. So the singular sets of linearly independent elements of an invariant \ccs over $\tor^2$ do not intersect. As before, the singular set of an element $\omega$ is non-empty and the number of its components is even, because the sign of the corresponding function $\varphi$ is different in adjacent regions of $\tor^2 \setminus \Sigma_\omega$ and because the components of $\Sigma_\omega$ do not bound disks, as there is an isotopy which carries any such component into itself, filling the whole torus on the way.
\finlem

\section{Invariant \cc s in dimension 3} \label{cc dim 3}
We are now looking for principal $\sph^1$-bundles over orientable, connected closed surfaces which admit invariant \cc s. It is in fact the base manifold that carries all the information with respect to this question. We have the following theorem\dop

\begin{theorem}\label{base-cercle} ---
Let $M$ be a principal circle-bundle over an orientable, connected closed surface $B$. There exists an invariant \ccs on $M$ if and only if the base space $B$ is either the 2-sphere or the 2-torus.
\end{theorem}

{\bf Proof (necessary part)\dop} Here we prove that only principal $\sph^1$-bundles over $\sph^2$ or $\tor^2$ might carry an invariant \cc. The converse will be proved in a constructive way in paragraph \ref{constructions}.

Let $\omega_1$ and $\omega_2$ generate an invariant \ccs on $M$. The singular sets $\Sigma_{\omega_1}$ and $\Sigma_{\omega_2}$ are unions of topological circles on $B$. By Corollary \ref{S1 fibr}, the complement of their intersection $\Sigma$ fibres over $\sph^1$.

If $\Sigma$ is empty, $B$ is an orientable, connected closed surface which fibres over $\sph^1$, thus the torus $\tor^2$.

If $\Sigma$ is non-empty, it is a finite union of points, as we have seen in paragraph \ref{ens sing}. According to the Fibration theorem \ref{thm fibr}, each component of $\Sigma$ has a neighborhood which is diffeomorphic to a disk and each radius of this disk corresponds to a different fibre of the associated knotted fibration over $\sph^1$. As each fibre starts and ends at two different components of $\Sigma$, $\Sigma$ has at least two components. $B \setminus \Sigma$ is a locally trivial fibration over $\sph^1$, so for reasons of continuity all fibres starting at a given component of $\Sigma$ end at the same component of $\Sigma$. With these restrictions the fibred surface can only be a twice punctured sphere, that is, $B=\sph^2$.
\finthm

\section{Invariant \cs s in dimension 3}
For invariant {\it \cs s} on principal circle-bundles over orientable, connected closed surfaces, the possibilities are even more restricted, as any \css contains \cc s. In fact, \cs s do not exist on principal circle-bundles over 2-tori, only on those over 2-spheres. We have the following theorem, analogous to the one in section \ref{cc dim 3}\dop

\begin{theorem}\label{base-sphere} ---
Let $M$ be a principal circle-bundle over a connected, orientable closed surface $B$. There exists an invariant \css on $M$ if and only if the base space $B$ is the 2-sphere.
\end{theorem}

Here again, we will only prove that the condition is necessary, whereas the constructions in paragraph \ref{constructions} will prove that it is also sufficient.

\vspace{0.4cm}
{\bf Proof of theorem \ref{base-sphere} (necessary part)\dop} By Theorem \ref{base-cercle} it is clear that the base space can only be $\sph^2$ or $\tor^2$. Let us see why it can not be the torus.

Assume that $\omega_1,\,\omega_2$ and $\omega_3$ generate an invariant \css on a principal circle bundle $M$ over $\tor^2$, with $\omega_i=\pi^*(\eta_i)+\pi^*(\varphi_i)\,\alpha,$ for $i=1,2,3$. We show that there exist two forms in the \css whose singular sets intersect, which contradicts Proposition \ref{descr cc}.

Remember that the singular set of any element of the \css is a non-empty union of topological circles on $B$, as we have seen in paragraph \ref{ens sing}.

Let $\Gamma$ be one connected component of $\Sigma_{\omega_1}$. By Proposition \ref{descr cc}, the singular sets of two elements of the \css do not intersect, so the function $\varphi$ associated to an arbitrary element of the \css is everywhere non-zero on $\Gamma$, except for $\varphi_1$ and $-\varphi_1$. Let $\omega$ be an element of $\mathcal{S}^2_c\{\omega_1,\,\omega_2,\,\omega_3\}$ other than $\omega_1$ and $-\omega_1$ and let
\begin{equation} \label{phi}
\varphi=\slas\,\varphi_i, \quad\textrm{where } \sla=1 \textrm{ and } \lambda_1\notin \{-1,1\},
\end{equation} 
be the associated function on the base space. We can assume that $\varphi$ is positive on $\Gamma$. Then there is a path on the unit sphere of triples $(\lambda_1,\,\lambda_2,\,\lambda_3)$ relating the coefficients of $\varphi$ in expression (\ref{phi}) to those of $-\varphi$ which does not take the values $(1,0,0)$ or $(-1,0,0)$, the coefficients of $\pm \varphi_1$. This path corresponds to a path in $\mathcal{S}^2_c\{\omega_1,\,\omega_2,\,\omega_3\}$ connecting $\omega$ and $-\omega$ and avoiding $\omega_1$ and $-\omega_1$. Thus we continuously transform $\varphi$, which is positive on $\Gamma$, into $-\varphi$, which is negative on $\Gamma$. Then there is an intermediate function $\tilde{\varphi}$, which corresponds to some element $\tilde{\omega}$ of $\mathcal{S}^2_c\{\omega_1,\,\omega_2,\,\omega_3\}$ and which has zeros on $\Gamma$. Thus $\Gamma$ and $\Sigma_{\tilde{\omega}}$ intersect.
\finthm

\section{Construction of examples} 
\label{constructions}
In the preceding paragraphs we proved necessary conditions for the existence of invariant \cc s and \cs s on principal $\sph^1$-bundles over surfaces. We will now see that they are also sufficient. On any manifold mentioned in Theorems \ref{base-cercle} and \ref{base-sphere}, that is, on principal $\sph^1$-bundles over $\sph^2$ or $\tor^2$, examples can be constructed.

Consider a principal circle-bundle $M$ over $\sph^2$ with bundle map $\pi$ and connection form $\alpha$ and three functions $\varphi_1,\,\varphi_2$ and $\varphi_3$ on $\sph^2$ such that $\sum \varphi_i^2=1$ everywhere. We define the following 1-forms on $M$\dop

\begin{eqnarray*}
\left\{ \begin{array}{rcl}
\omega_1&=&\pi^*(\varphi_2\,d\varphi_3-\varphi_3\,d\varphi_2)+\pi^*(\varphi_1)\,k\,\alpha\\
\omega_2&=&\pi^*(\varphi_3\,d\varphi_1-\varphi_1\,d\varphi_3)+\pi^*(\varphi_2)\,k\,\alpha\\
\omega_3&=&\pi^*(\varphi_1\,d\varphi_2-\varphi_2\,d\varphi_1)+\pi^*(\varphi_3)\,k\,\alpha,
\end{array} \right.
\end{eqnarray*}
where $k$ is a real number.

\begin{lemma} \label{formule magique} ---
Consider the form $\omega=\lambda_1\,\omega_1+\lambda_2\,\omega_2+\lambda_3\,\omega_3$ for some $\lambda_1,\lambda_2,\lambda_3$ such that $\sum \lambda_i^2=1$. Then we have the formula\dop
{\setlength\arraycolsep{2pt}
\begin{eqnarray*} \omega\wedge d\omega&=&k\,\pi^*((\sum_{i=1}^3\lambda_i\,\varphi_i)(\lambda_1\,d\varphi_2\wedge d\varphi_3+\lambda_2\,d\varphi_3\wedge d\varphi_1+\lambda_3\,d\varphi_1\wedge d\varphi_2))\wedge\alpha\\
&&+k\,\pi^*(\Phi^*(\Omega))\wedge\alpha+k^2\,\pi^*((\sum_{i=1}^3\lambda_i\,\varphi_i)^2)\alpha\wedge d\alpha
\end{eqnarray*}}
where $\Omega$ is the Euclidean volume form on $\sph^2$ and $\Phi=(\varphi_1,\varphi_2,\varphi_3): \sph^2 \rightarrow \sph^2$.
\end{lemma}

{\bf Proof\dop} This is an easy, though lengthy calculation using the identities $\sum\lambda_i^2=1$, $\sum \varphi_i^2=1$ and $\Phi^*(\Omega)=\varphi_3\,d\varphi_1\wedge\,d\varphi_2+\varphi_1 d\varphi_2\wedge\,d\varphi_3+\varphi_2 d\varphi_3\wedge d\varphi_1$ and the fact that $d\alpha$ is an horizontal form. The only form involved which is not a pullback from $\sph^2$ is $\alpha$.
\finlem

This lemma gives us the possibility to construct contact spheres on principal $\sph^1$-bundles over $\sph^2$. Indeed, up to an adjustment of the connection form $\alpha$, the forms $(\omega_1,\,\omega_2,\,\omega_3)$ given above generate a \css for the right choice of $k$ and of the functions $(\varphi_1,\,\varphi_2,\,\varphi_3)$.

Consider a function triple $(\varphi_1,\,\varphi_2,\,\varphi_3)$ such that $\sum_{i=1}^3 \varphi_i^2=1$ and which defines a knotted fibration of $\sph^2$ over $\sph^2$, i.e. the corresponding common singular set $\bigcap_{i=1}^3\varphi_i^{-1}(0)$ is empty and $\Phi=(\varphi_1,\,\varphi_2,\,\varphi_3)$ is of rank two everywhere. As the structural group of the fibration is abelian, the form $d\alpha$ is horizontal, thus there is a 2-form $\nu$ on the base space $B$, such that $d\alpha=\pi^*(\nu)$. Since $\Phi^*(\Omega)$ is a volume form on $\sph^2$, there is a function $f$ such that $\nu=f\,\Phi^*(\Omega)$ and some functions $C_1, C_2, C_3$, such that $d\varphi_2\wedge d\varphi_3=C_1\,\Phi^*(\Omega),\,d\varphi_3\wedge d\varphi_1=C_2\,\Phi^*(\Omega)$ and $d\varphi_1\wedge d\varphi_2=C_3\,\Phi^*(\Omega)$. The right hand expression in Lemma 2 now becomes
\begin{eqnarray} \label{f.m.2} k\big(\pi^*\big(\big(1+k\,f\,(\sum_{i=1}^3\lambda_i\,\varphi_i)^2+(\sum_{i=1}^3\lambda_i\,\varphi_i)(\lambda_1\,C_1+\lambda_2\,C_2+\lambda_3\,C_3)\big)\Phi^*(\Omega)\big)\wedge \alpha\big).
\end{eqnarray}

To make this form a volume form, we observe that the function $(\sum_{i=1}^3\lambda_i\,\varphi_i)(\lambda_1\,C_1+\lambda_2\,C_2+\lambda_3\,C_3)$ is bounded on $\sph^2$ and that it is thus sufficient to choose $k$ large enough and of the right sign, if $f$ is a non-vanishing function of constant sign. Next, we notice that adding the function $h=\frac{\int_{\sph^2}f\,\Phi^*(\Omega)}{\int_{\sph^2}\Phi^*(\Omega)}-f$ to $f$ makes $f$ a constant function. This is possible, for the following reason\dop the form $\eta=h\,\Phi^*(\Omega)$ satisfies $\int_{\sph^2}\eta=0$, so $\pi^*(\eta)$ is an exact form on $M$ and there is a 1-form $\xi$ such that $d\xi =\pi^*(\eta)$. The fibration of $M$ over $\sph^2$ is characterized by the value of $\int_{\sph^2}\nu$, where $\pi^*(\nu)$ is the differential of the connection form, so $\alpha+\xi$ is another connection form of the same fibration. Thus replacing $\alpha$ by $\alpha+\xi$ changes $f$ into a constant function, without modifying the fibration. 

Now assume $f$ to be constant. If $f$ is not identically zero, choosing $k$ large enough and of the same sign as $f$ makes the expression (\ref{f.m.2}) a volume form. If $f\equiv0$, then the manifold $M$ is equivariantly diffeomorphic to $\sph^2\times\sph^1$ and there are known examples of invariant \cs s on such manifolds (see \cite{GG2}, p.274). This gives us examples of invariant \cs s on all principal $\sph^1$-bundles over $\sph^2$.

\vspace{0.4cm} 
Let us now construct examples of invariant \cc s on principal $\sph^1$-bundles over $\tor^2$. In this case we consider the following 1-forms on $M$\dop

\begin{eqnarray*}
\left\{ \begin{array}{rcl}
\omega_1&=&\pi^*(\cos\theta_1\,d\theta_2)+\pi^*(\sin\theta_1)\,k\,\alpha\\
\omega_2&=&\pi^*(-\sin\theta_1\,d\theta_2)+\pi^*(\cos\theta_1)\,k\,\alpha,
\end{array} \right.
\end{eqnarray*}
$\theta_1$ and $\theta_2$ being pseudo-coordinates on $\tor^2$,$\,\alpha$ a connection form for the fibration of $M$ over $\tor^2$ and $k$ a real number. Setting $\omega=\lambda_1\,\omega_1+\lambda_2\,\omega_2$, we get the formula
\begin{eqnarray} \label{f.m.tore}
\omega\wedge d\omega&=&-k\,d\theta_1\wedge d\theta_2\wedge\alpha+k^2(\lambda_1\,\sin\theta_1+\lambda_2\,\cos\theta_1)^2\alpha\wedge d\alpha.
\end{eqnarray}
Similarly to the previously discussed case, there is a function $f$ such that $d\alpha=f\,d\theta_1\wedge d\theta_2$, and we can assume that $f$ is constant, up to a modification of $\alpha$ by addition of a 1-form whose differential comes from a basic 2-form integrating to zero over the torus. Thus, (\ref{f.m.tore}) becomes
\begin{eqnarray}
\omega\wedge d\omega&=&k(k\,f\,(\lambda_1\,\sin\theta_1+\lambda_2\,\cos\theta_1)^2-2)\,d\theta_1\wedge d\theta_2\wedge\alpha.
\end{eqnarray}
It is now obvious that it is enough to choose $k$ non-zero and such that $k\cdot f$ is not positive to make $\omega\wedge d\omega$ a volume form, for any coefficients $(\lambda_1,\,\lambda_2)$.

It may be interesting to see that a \ccs obtained by this construction is taut only if $f$ is identically zero, that is, for $M=\tor^3$.
\vspace{0.4cm}

The above construction leads to the following proposition, which gives examples of invariant \cc s and \cs s on all manifolds mentioned in Theorems \ref{base-cercle} and \ref{base-sphere}. Thus it finishes the proves of these two theorems.

\begin{prop} \label{existence} ---
Let $M$ be a circle-bundle over $\sph^2$ with connection form $\alpha$ and $\varphi_1,\,\varphi_2$ and $\varphi_3$ three functions on $\sph^2$ such that $\sum \varphi_i^2=1$ everywhere and which define a knotted fibration of $\sph^2$ over $\sph^2$. Then there exist a real number $k$ and a 1-form $\xi$ on $M$, such that the 1-forms
\begin{eqnarray*}
\left\{ \begin{array}{rcl}
\omega_1&=&\pi^*(\varphi_2\,d\varphi_3-\varphi_3\,d\varphi_2)+\pi^*(\varphi_1)\,k(\alpha+\xi)\\
\omega_2&=&\pi^*(\varphi_3\,d\varphi_1-\varphi_1\,d\varphi_3)+\pi^*(\varphi_2)\,k(\alpha+\xi)\\
\omega_3&=&\pi^*(\varphi_1\,d\varphi_2-\varphi_2\,d\varphi_1)+\pi^*(\varphi_3)\,k(\alpha+\xi)
\end{array} \right.
\end{eqnarray*}
generate an invariant \css on $M$.

Similarly, if $M$ is a circle-bundle over $\tor^2$ with connection form $\tilde{\alpha}$ and pseudo-coordinates $\theta_1,\,\theta_2$ on $\tor^2$, then there exist a real number $\tilde{k}$ and a 1-form $\tilde{\xi}$ on $M$, such that the 1-forms
\begin{eqnarray*}
\left\{ \begin{array}{rcl}
\omega_1&=&\pi^*(\cos\theta_1\,d\theta_2)+\pi^*(\sin\theta_1)\,\tilde{k}\,(\tilde{\alpha}+\tilde{\xi})\\
\omega_2&=&\pi^*(-\sin\theta_1\,d\theta_2)+\pi^*(\cos\theta_1)\,\tilde{k}\,(\tilde{\alpha}+\tilde{\xi})
\end{array} \right.
\end{eqnarray*}
generate an invariant \ccs on $M$.
\end{prop}

\section{Non-existence of contact $p$-spheres in dimension $4n+1$}

We now consider the situation in higher dimensions. It quickly becomes obvious that in dimension 5 it is not easy to construct \cc s, although there are, of course, 5-manifolds admitting contact forms. In fact, we prove the following\dop

\begin{theorem} \label{thm dim 5} ---
On 5-dimensional manifolds, and more generally on $(4n+1)$-dimensional manifolds, \cc s and a fortiori contact $p$-spheres do not exist for $p\geq 1$.
\end{theorem}

{\bf Proof\dop} Let $(\omega_1,\,\omega_2)$ be a pair of contact forms on a 5-dimensional manifold $M$. It generates a \ccs if $\omega=\lambda_1\,\omega_1+\lambda_2\,\omega_2$ is a contact form for any $(\lambda_1,\,\lambda_2)$ with $\lambda_1^2+\lambda_2^2=1$, that is, if
$$
\omega\wedge\,(d\omega)^2=\sum_{i,j,\,k=1}^2\lambda_i\,\lambda_j\,\lambda_k\,(\omega_i\wedge\,d\omega_j\wedge\,d\omega_k)
$$
is everywhere non-zero.

Now fix a point $x\in M$ and five linearly independent tangent vectors $$(v_1,\ldots,\,v_5)\in {\cal T}_x(M)$$ and consider the function 
\begin{eqnarray*}
\R^2&\longrightarrow&\R\\
(\lambda_1,\,\lambda_2)&\longmapsto&\sum_{i,j,k=1}^2\lambda_i\,\lambda_j\,\lambda_k\,(\omega_i\wedge\,d\omega_j\wedge\,d\omega_k)_x(v_1,\ldots,\,v_5),
\end{eqnarray*}
which is a homogeneous polynomial function of degree 3. It has zeros on the unit circle of $\R^2$ (if it is positive at some point of the circle, it is negative at its antipode), so $\mathcal{S}^1_c\{\omega_1,\,\omega_2\}$ cannot be a \cc, as for the corresponding coefficients $(\lambda_1,\lambda_2)$, $\omega\wedge\,(d\omega)^2$ is not a volume form.

In general, pairs of contact forms in dimension $4n+1$ give us polynomial functions of degree $2n+1$, which is odd, so they necessarily vanish on the unit circle of $\R^2$. Thus, \cc s do not exist in these dimensions.

On the other hand, in dimension $4n-1$, these polynomial functions are of degree $2n$, which is even, so there is no restriction to the existence of \cc s in these dimensions.
\finthm

\section{Round contact $p$-spheres} \label{round cs} 

Let us now consider \cps s with a special regularity property. H. Geiges and J. Gonzalo call a \cpss \emph{taut} if all elements yield the same volume form. We introduce a more geometrical property which is equivalent to tautness in dimension 3, as we shall see.

\subsection{Definition and examples}

\begin{defi} ---
Let $\{\omega_1,\dots,\omega_{p+1}\}$ generate a contact $p$-sphere and let $R_1,\dots,R_{p+1}$ be the corresponding Reeb vector fields. Then   $\mathcal{S}^p_c\{\omega_1,\dots,\omega_{p+1}\}$ is called a \emph{round} contact $p$-sphere if any form $\sum_{i=1}^{p+1} \lambda_i\,\omega_i$ with $\sum_{i=1}^{p+1} \lambda_i^2=1$ admits $\sum_{i=1}^{p+1} \lambda_iR_i$ as its Reeb vector field.
\end{defi}

So the Reeb vector field of a linear combination of elements of a round \cpss is the corresponding linear combination of the Reeb vector fields of these elements. This is a strong restriction and quite useful for computations.

We have the following characterization of round contact $p$-spheres\dop

\begin{lemma}\label{CNS rond} ---
Let $\{\omega_1,\dots,\omega_{p+1}\}$ generate a contact $p$-sphere on some manifold $M$ and let $R_1,\dots,R_{p+1}$ be the corresponding Reeb vector fields. Then $\mathcal{S}^p_c\{\omega_1,\dots,\omega_{p+1}\}$ is round if and only if the following conditions are satisfied\dop
\begin{eqnarray*} \begin{array}{lcl}
i)\quad \ \, \omega_i(R_j)+\omega_j(R_i)=0,\textrm{ for }i,j\in\{1,\dots,p+1\},\,i \neq j\\
ii)\quad R_i \lrcorner\,d\omega_j+R_j\lrcorner\,d\omega_i=0,\textrm{ for }i,j\in\{1,\dots,p+1\}.
\end{array} \end{eqnarray*}
\end{lemma}

{\bf Proof\dop} Consider a contact form $\omega=\sum_{i=1}^{p+1}\lambda_i\,\omega_i$ with normalized coefficients. $R=\sum_{i=1}^{p+1}\lambda_i\,R_i$ is the corresponding Reeb vector field if and only if 
\begin{eqnarray*} \begin{array}{lcl}
a)\quad (\sum_{i=1}^{p+1}\lambda_i\,\omega_i)(\sum_{i=1}^{p+1}\lambda_i\,R_i)=1\\[1mm]
b)\quad (\sum_{i=1}^{p+1}\lambda_i\,R_i) \lrcorner (\sum_{i=1}^{p+1}\lambda_i\,d\omega_i)=0
\end{array} \end{eqnarray*}
Condition $a)$ means that $\sum_{i<j} \lambda_i\lambda_j(\omega_i(R_j)+\omega_j(R_i))=0$. If $\mathcal{S}^p_c\{\omega_1,\dots,\omega_{p+1}\}$ is a round contact $p$-sphere, this equality is true for any normalized $(p+1)$-tuple of coefficients $(\lambda_1,\dots,\lambda_{p+1})$. Considering the left-hand side as a homogeneous polynomial on $\sph^p$, we conclude that all coefficients must be zero, that is, condition $i)$ is satisfied. The converse is of course true.

The equivalence of conditions $b)$ and $ii)$ is proved in the same way.
\finlem

{\bf Remark\dop}This lemma shows that a \cpss is round if and only if every \ccs it contains is round. The same is true for tautness. So to show that a \cpss is round or taut, it is enough to examine \cc s for this property.
\vspace{0.4cm}

{\bf Examples\dop}
\begin{enumerate}
\item We have mentioned an example of a round contact $2$-sphere on $\sph^3$ in paragraph \ref{exemples}. Indeed, consider $\sph^3$ as the unit sphere of the group of quaternions $\ha$ and the \css generated by $\omega_1,\,\omega_2$ and $\omega_3$ induced on $\sph^3$ by
\begin{eqnarray*}
\left\{ \begin{array}{rcl}
\tilde{\omega}_1&=&<iq,dq>\\
\tilde{\omega}_2&=&<jq,dq>\\
\tilde{\omega}_3&=&<kq,dq>.
\end{array} \right.
\end{eqnarray*}
The corresponding Reeb vector fields are
\begin{eqnarray*}
\left \{ \begin{array}{rcl}
R_1&=&iq\\
R_2&=&jq\\
R_3&=&kq.
\end{array} \right.
\end{eqnarray*}

Using Lemma \ref{CNS rond}, we can check that this \css is indeed round\dop

In this example, the Reeb vector fields of $\omega_1,\,\omega_2$ and $\omega_3$ are also their dual vector fields, so $\omega_i(R_j)=\delta_{ij}$. So condition $i)$ is satisfied.

On the other hand, as $d \tilde{\omega}_1=-d\bar{q}\wedge idq$ and $d \tilde{\omega}_2=-d\bar{q}\wedge jdq$, we have 
\begin{eqnarray*}
R_2\lrcorner\,d\tilde{\omega}_1&=&\bar{q}\,j\,i\,dq+d\bar{q}\,i\,j\,q\\
R_1\lrcorner\,d\tilde{\omega}_2&=&\bar{q}\,i\,j\,dq+d\bar{q}\,j\,i\,q\\
&=&-(\bar{q}\,j\,i\,dq+d\bar{q}\,i\,j\,q).
\end{eqnarray*}
So $R_2\lrcorner\,d\omega_1+R_1\lrcorner\,d\omega_2=0$ and the remaining relations of condition $ii)$ are obtained in the same way.

\item There is also an example of a round \ccs on $\tor^3$. Indeed, the forms

\begin{eqnarray*}
\left\{ \begin{array}{rcl}
\omega_1&=&\cos\theta_1\,d\theta_2+\sin\theta_1\,d\theta_3\\
\omega_2&=&-\sin\theta_1\,d\theta_2+\cos\theta_1\,d\theta_3
\end{array} \right.
\end{eqnarray*}

generate a \ccs and their Reeb vector fields are respectively
\begin{eqnarray*}
\left\{ \begin{array}{rcl}
R_1&=&\cos\theta_1\,\frac{\partial}{\partial\theta_2}+\sin\theta_1\,\frac{\partial}{\partial\theta_3}\\
R_2&=&-\sin\theta_1\,\frac{\partial}{\partial\theta_2}+\cos\theta_1\,\frac{\partial}{\partial\theta_3}.
\end{array} \right.
\end{eqnarray*}

Once again, we have $\omega_i(R_j)=\delta_{ij}$, for $i,\,j=1,2$. So condition $i)$ of Lemma \ref{CNS rond} is satisfied. Moreover, we have
\begin{eqnarray*}
\left \{ \begin{array}{rcl}
d\omega_1&=&-\sin\theta_1\,d\theta_1\wedge d\theta_2+\cos\theta_1\,d\theta_1\wedge d\theta_3\\
d\omega_2&=&-\cos\theta_1\,d\theta_1\wedge d\theta_2-\sin\theta_1\,d\theta_1\wedge d\theta_3,
\end{array} \right.
\end{eqnarray*}
so $R_1\lrcorner d\omega_2=d\theta_1=-R_2\lrcorner d\omega_1$ and condition $ii)$ is also satisfied.
\end{enumerate}

\subsection{Round $\Leftrightarrow$ taut in dimension 3} 
We now prove that in dimension 3, roundness is indeed equivalent to tautness. This gives us the possibility to have different viewpoints on the same property. Let us first prove the following lemma:

\begin{lemma} \label{Reeb ind} ---
The Reeb vector fields of two generating elements of a \ccs are everywhere linearly independent.
\end{lemma}

{\bf Proof\dop} Let $\omega_1$ and $\omega_2$ generate a \ccs and let $R_1$ and $R_2$ be the corresponding Reeb vector fields. Assume that in some point $p$, $(R_1)_p$ and $(R_2)_p$ are parallel. Then $d\omega_1(R_1)_p=0$ and $d\omega_2(R_1)_p=0$, so $d\omega(R_1)_p=0$, for any linear combination $\omega$ of $\omega_1$ and $\omega_2$. Thus the Reeb vector fields of all elements of the \ccs are parallel to $R_1$. Now, the Reeb vector field of $-\omega_1$, which is an element of the \cc, is $-R_1$, and as the Reeb vector field depends continuously on the coefficients, there must be some form in the \ccs with a Reeb vector field of length zero. This is of course impossible.
\finlem

\begin{theorem} ---
On a 3-dimensional manifold $M$, a \ccs (resp. a \cs) is taut if and only if it is round.
\end{theorem} 

{\bf Proof\dop} Let us first see that taut \cc s are round. Let $\mathcal{S}^1_c\{\omega_1,\,\omega_2\}$ be a taut \cc, that is, satisfying
\begin{eqnarray} \label{taut}
\left \{ \begin{array}{lll}
\omega_1 \wedge d\omega_1&=&\omega_2 \wedge d\omega_2 \\
\omega_1 \wedge d\omega_2&=&-\omega_2 \wedge d\omega_1
\end{array} \right.
\end{eqnarray}
and let $R_1$ and $R_2$ be the corresponding Reeb vector fields. Applying the first equation of (\ref{taut}) to the couple of vectors $(R_1,\,R_2)$, we get $$R_1 \lrcorner \,  d\omega_2+R_2 \lrcorner \, d\omega_1=0,$$ which is the second condition of roundness of Lemma \ref{CNS rond}. Applying the second equation of (\ref{taut}) to the couple $(R_1,\,R_2)$ and using the previous relation, we get $$ (R_1 \lrcorner \, d\omega_2)(\omega_1(R_2)+\omega_2(R_1))=0,$$ which gives us the first condition of roundness of Lemma \ref{CNS rond}, if we are sure that $R_1 \lrcorner \, d\omega_2$ never vanishes. But this is granted by Lemma \ref{Reeb ind}. So taut \cc s are round.

\vspace{0.4cm}
It is now immediate that taut \cs s are round, because the roundness conditions of Lemma \ref{CNS rond} carry on pairs of generating contact forms. So if two generating forms of a \css satisfy the tautness condition (\ref{taut}), then they also satisfy the roundness condition, as we have just seen.

\vspace{0.4cm}
Let us now see that round \cs s are taut. Let $(\omega_1,\,\omega_2,\,\omega_3)$ generate a round \css and consider an element $\omega=\sum_{i=1}^3\lambda_i\,\omega_i$. By assumption, the corresponding Reeb vector field is $R=\sum_{i=1}^3\lambda_i\,R_i$, where $R_i$ is the Reeb vector field of the form $\omega_i$. Thus non-trivial linear combinations of $R_1,\,R_2$ and $R_3$ never vanish, so $R_1,\,R_2$ and $R_3$ are everywhere linearly independent vector fields, hence a coordinate system on $M$. To prove that $\mathcal{S}^2_c\{\omega_1,\,\omega_2,\,\omega_3\}$ is taut, it is enough to show that $\omega\wedge d\omega (R_1,\,R_2,\,R_3)$ is independent of the coefficients $(\lambda_1,\,\lambda_2,\,\lambda_3)$. An evaluation of the relation $R_2 \lrcorner\,d\omega_1+R_1\lrcorner\,d\omega_2=0$ on the vector field $R_3$ yields
$$d\omega_1(R_2,\,R_3)=d\omega_2(R_3,\,R_1),$$
and by an analogous evaluation we get 
$$d\omega_1(R_2,\,R_3)=d\omega_3(R_1,\,R_2).$$

So we have
\begin{eqnarray*}
\omega\wedge d\omega (R_1,\,R_2,\,R_3)&=&(\lambda_1\,\omega_1+\lambda_2\,\omega_2+\lambda_3\,\omega_3)\wedge(\lambda_1\,d\omega_1+\lambda_2\,d\omega_2+\lambda_3\,d\omega_3)(R_1,\,R_2,\,R_3)\\
&=&\lambda_1^2d\omega_1(R_2,\,R_3)+\lambda_2^2d\omega_2(R_3,\,R_1)+\lambda_3^2d\omega_3(R_1,\,R_2)\\
&=&d\omega_1(R_2,\,R_3),
\end{eqnarray*}
where the second step is due to the preceding relations and Lemma \ref{CNS rond}. This expression does not depend on the coefficients $(\lambda_1,\,\lambda_2,\,\lambda_3)$.

\vspace{0.4cm}
It remains to show that round \cc s are taut in dimension $3$. So let $\omega_1$ and $\omega_2$ generate a round \ccs and let $R_1$ and $R_2$ be the corresponding Reeb vector fields. By Lemma \ref{Reeb ind}, $R_1$ and $R_2$ are everywhere linearly independent, so we can find a third vector field $X$ to make $(R_1,\,R_2,\,X)$ a basis of the tangent bundle. An evaluation of the relation \mbox{$R_1 \lrcorner \, d\omega_2+R_2 \lrcorner \, d\omega_1=0$} on $X$ gives us $$d\omega_1(R_2,\,X)=d\omega_2(X,\,R_1).$$ Now an analogous calculation as in the case of \cs s shows that the \ccs is taut.
\finthm

\subsection{Round vs taut in higher dimensions}

In dimensions higher than 3, roundness and tautness are not equivalent. In fact, from dimension 7 on, there is an essential difference between the two notions, due to the powers of the differentials which appear in the definition of tautness and not in the definition of roundness. To illustrate this difference, let us consider the situation in dimension 7.

A \ccs $\Sc^1\{\omega_1,\,\omega_2\}$ defined on a $7$-dimensional manifold is taut if the following four equalities are satisfied\dop
\begin{eqnarray}
\left \{ \begin{array}{lcl} \label{taut dim 7}
\omega_1 \wedge (d\omega_1)^3-\omega_2 \wedge (d\omega_2)^3 &=&0 \\
3\,\omega_2 \wedge d\omega_2 \wedge (d\omega_1)^2+3\,\omega_1 \wedge d\omega_1 \wedge (d\omega_2)^2-2\,\omega_1 \wedge (d\omega_1)^3 &=&0 \\
\omega_1 \wedge (d\omega_2)^3+3\,\omega_2 \wedge d\omega_1 \wedge (d\omega_2)^2 &=&0 \\
\omega_2 \wedge (d\omega_1)^3+3\,\omega_1 \wedge d\omega_2 \wedge (d\omega_1)^2 &=&0.
\end{array} \right .
\end{eqnarray}

On the other hand, if $\Sc^1\{\omega_1,\,\omega_2\}$ is taut and round, we have the following equalities, which are necessary, but not sufficient\dop
\begin{eqnarray}
\left \{ \begin{array}{lcl} \label{rond+taut}
(d\omega_1)^3-3\,d\omega_1 \wedge (d\omega_2)^2 &=& 0 \\
(d\omega_2)^3-3\,d\omega_2 \wedge (d\omega_1)^2 &=& 0.
\end{array} \right .
\end{eqnarray}
This can be seen as follows\dop
Let $\omega_\theta=\cos\theta\,\omega_1+\sin\theta\,\omega_2$ be an element of $\Sc^1\{\omega_1,\,\omega_2\}$ and note $\Omega=\omega_1 \wedge (d\omega_1)^3.$ As $\Sc^1\{\omega_1,\,\omega_2\}$ is taut, we have $\omega_\theta \wedge (d\omega_\theta)^3=\Omega,$ and thus $(d\omega_\theta)^3=R_\theta \lrcorner \Omega,$ where $R_\theta$ is the Reeb vector field of $\omega_\theta$. This relation leads to the conditions (\ref{rond+taut}).

\vspace{0.4cm}
The systems (\ref{taut dim 7}) et (\ref{rond+taut}) are not of the same nature, as the equations of the first one are of degree 7 and those of the second one are of degree 6. So in general, a \ccs which satisfies (\ref{taut dim 7}) will not satisfy (\ref{rond+taut}). In dimension 3, the corresponding equations of (\ref{rond+taut}) are trivial, this is why we have equivalence of tautness and roundness in dimension 3.

\vspace{0.4cm}
The following counter-examples prove that in dimension 7 there is no implication between tautness and roundness\dop
\begin{enumerate}
\item The \ccs on $\R^7$ generated by the forms 
\begin{eqnarray*}
\left \{ \begin{array}{lcl}
\omega_1 &=& x_1\,dx_2+x_3\,dx_4+x_5\,dx_6+dx_7 \\
\omega_2 &=& -(x_5+x_6)\,dx_3-x_5\,dx_4+(x_1+x_3)\,dx_6+x_1\,dx_7-dx_2
\end{array} \right .
\end{eqnarray*}
is round and not taut.

\item The \ccs on $\R^7$ generated by the forms 
\begin{eqnarray*}
\left \{ \begin{array}{lcl}
\omega_1 &=& x_1\,dx_2+x_3\,dx_4+x_5\,dx_6+dx_7 \\
\omega_2 &=& x_5\,dx_4-x_3\,dx_6+(x_1+x_3)\,dx_7-dx_2
\end{array} \right .
\end{eqnarray*}
is taut and not round.
\end{enumerate}

\section{Examples of contact $p$-spheres in higher dimensions}

We have seen that \cc s do not exist on manifolds of dimension $4n+1$. In the other odd dimensions, however, many interesting examples of contact $p$-spheres can be found. A first family of manifolds where contact $p$-spheres are likely to be found easily are the spheres of dimension $4n-1$.

On the spheres, we have a natural upper bound of the size of contact $p$-spheres, given by Adams' formula. According to Adams, there do not exist more than $\rho(n)$ continuous unitary vector fields on $\sph^{n-1}$, which are everywhere linearly independent, where $$\rho(n)=2^c+8d-1,\textrm{ with } n=odd \cdot 2^{c+4d}, \; c\leq 3.$$
This means of course that there can not be more than $\rho(n)$ everywhere linearly independent contact forms neither. So on $\sph^{4n-1}$, contact $p$-spheres can only exist for $p \leq (\rho(4n)-1)$.

On the other hand, there are works of B. Eckmann, relying on ideas of A. Hurwitz and J. Radon, leading to the following theorem\dop

\begin{theorem} ---
On $\sph^{4n-1}$, there exists a contact $(\rho(4n)-1)$-sphere, for $n \geq 1$.
\end{theorem}

{\bf Proof\dop} According to B. Eckmann (see \cite{Eck}), there exist, for each integer $n$, $\rho(4n)$ antisymmetric matrices ${A_1,\dots,A_{\rho(4n)}}$ of $\mathcal{O}(4n,\R)$, such that the vectors ${(A_1x)_x,\dots,(A_{\rho(4n)}x)_x}$ are linearly independent everywhere on the unit sphere of $\R^{4n}$, that is, Eckmann gives a realization of the maximum number of such vector fields given by Adam´s formula. The same matrices can also be used to construct a contact $(\rho(4n)-1)$-sphere on $\sph^{4n-1}$, as we will see now.

The matrices ${A_1,\dots,A_{\rho(4n)}}$ satisfy the relations
\begin{eqnarray}
\left \{ \begin{array}{lcr}
A_iA_j+A_jA_i=0, \; i \neq j \\
A_i^2=-Id, \; i=1, \dots ,\rho(4n),
\end{array} \right .
\end{eqnarray}
that is, they generate a Clifford algebra.

Now define $\rho(4n)$ 1-forms on $\R^{4n}$ by $$(\tilde{\omega}_i)_x=<A_ix,dx>=\sum_{r,s=1}^{4n} a^i_{rs}x_s dx_r,$$
where $A_i=(a^i_{rs})_{r,s=1\dots 4n}$.

These forms induce some forms $\omega_i$ on $\sph^{4n-1}$, which are contact forms: as each matrix $A_i$ is orthogonal and antisymmetric, there is an orthogonal basis of $\R^{4n}$ in which $A_i$ becomes 
\begin{eqnarray} \label{matrice}
\left ( \begin{array}{rrrrrr}
	0&-1&&&&\\
	1&0&&&&\\
	&&.&&&\\
	&&&.&&\\
	&&&&0&-1\\
	&&&&1&0
\end{array} \right ) .
\end{eqnarray}

In this basis, we have $$<A_ix,dx>=\sum_{i=1}^{2n} (x_{2i-1} dx_{2i}-x_{2i} dx_{2i-1}),$$ which induces a contact form on $\sph^{4n-1}$.

Moreover, and quite surprisingly, these $\rho(4n)$ contact forms generate a $(\rho(4n)-1)$-\css on $\sph^{4n-1}$. 
Each matrix $$\sum_{i=1}^{\rho(4n)} \lambda_i A_i \textrm{ with } \sum_{i=1}^{\rho(4n)} \lambda_i^2=1$$ is indeed antisymmetric and it is orthogonal as well, as 
\begin{eqnarray*}
^t(\sum_{i=1}^{\rho(4n)} \lambda_i A_i)(\sum_{i=1}^{\rho(4n)} \lambda_i A_i)&=&-\sum_{i=1}^{\rho(4n)} \lambda_i^2 A_i^2 - \sum_{i<j} \lambda_i \lambda_j (A_iA_j+A_jA_i)\\
&=&Id.
\end{eqnarray*}

So the corresponding 1-form $<\sum_{i=1}^{\rho(4n)} \lambda_i A_i x,dx>$ induces a contact form, as we have already seen.
\finthm

This theorem gives us examples of contact $p$-spheres on $\sph^{4n-1}$ with $p \geq 2$ for any $n\geq 1$. According to Adams' formula, we have $\rho(4n)=2^c+8d-1$, where $4n=odd \cdot 2^{c+4d},\,c\leq 3$. In our situation, we have $c \geq 2$ or $d \geq 1$ and in both cases $\rho(4n)\geq 3$, which means that on $\sph^{4n-1}$, there is a contact $p$-sphere generated by at least $3$ contact forms.

It may be interesting to see that these examples are round and taut. 

Indeed, if we consider a \ccs $\Sc^1\{\omega_1,\,\omega_2\}$ which is contained in one of these examples on some sphere $\sph^{4n-1}$, the Reeb vector field of a form $\omega_i$ induced by $(\tilde{\omega}_i)_x=<A_ix,dx>$ on $\sph^{4n-1}$ is $(R_i)_x=(A_i\,x)_x$. Then we have
\begin{eqnarray*}
	\left \{ \begin{array}{lclclcl}
	\tilde{\omega}_1(R_2) &=& <A_1\,x,A_2\,x> &=& <-A_2\,A_1\,x,x> &=& 0 \\
	\tilde{\omega}_2(R_1) &=& <A_2\,x,A_1\,x> &=& <-A_1\,A_2\,x,x> &=& 0 \\
	(R_2 \lrcorner\, d\tilde{\omega}_1)_x &=& ^t x\,A_1\,A_2\,dx &=& -\,^t x\,A_2\,A_1\,dx &=& -(R_1 \lrcorner\, d\tilde{\omega}_2)_x,
	\end{array} \right .
\end{eqnarray*}
so $\Sc^1\{\omega_1,\,\omega_2\}$ is round.

Now consider an element of the \cc $$\omega=\lambda_1\,\omega_1+\lambda_2\,\omega_2 \quad \textrm{with} \quad \lambda_1^2+\lambda_2^2=1.$$ 
$\omega$ is induced by $(\tilde{\omega})_x=<(\lambda_1\,A_1+\lambda_2\,A_2)x,dx>$. There is an orthogonal, thus volume preserving coordinate change which transforms the matrix $\lambda_1\,A_1+\lambda_2\,A_2$ into the matrix (\ref{matrice}). Thus $\omega$ defines the same volume form as $\omega_1$ and $\Sc^1\{\omega_1,\,\omega_2\}$ is taut.

\section{Contact spheres and Sasakian 3-structures}

Another way to find examples of contact spheres on higher dimensional manifolds is to consider Sasakian geometry, and more precisely Sasakian 3-structures. As a Sasakian 3-structure is defined by three contact forms with special properties, we can ask if all non-trivial linear combinations of these forms are still contact forms, that is, if they define a contact sphere. In that case, we can study the regularity properties of such a contact sphere. Proposition \ref{3-structures}, which has been suggested to me by the referee, answers these questions. Let us first recall the definition of a Sasakian 3-structure (see \cite{Bl}).
\vspace{4mm}

Let $(M,g)$ be a Riemannian manifold which carries a contact form $\omega$ and let $R$ be the Reeb vector field associated to $\omega$. We define the tensor field $\varphi$ of type $(1,1)$ by
\begin{equation*} \label{tenseur phi} g(X,\varphi(Y))=\frac{1}{2} \, d\omega(X,Y). \end{equation*}
$\omega$ is a metric contact form if it satisfies
$$g(\varphi(X),\varphi(Y))=g(X,Y)-\omega(X)\,\omega(Y),$$
or equivalently
\begin{eqnarray*} \label{forme métrique}
\left\{ \begin{array}{rcl}
\varphi^2 &=& -I + \omega \otimes R \\
\omega(Y) &=& g(R,\,Y),\; \textrm{for any} \; Y.
\end{array} \right. 
\end{eqnarray*}
$(\varphi,\,R,\,\omega)$ defines a Sasakian structure if it satisfies
\begin{equation*} \label{Sasaki} \left[\varphi,\varphi \right] + d\omega \otimes R = 0, \end{equation*}
where $ \left[ \cdot, \cdot \right] $ is the Nijenhuis bracket, defined by
$$ \left[ T, T \right](X,Y) = T^2 \left[ X, Y \right] + \left[ TX, TY \right] - T \left[ TX, Y \right] - T \left[ X, TY \right].$$

Three given Sasakian structures $\omega_1,\,\omega_2$ et $\omega_3$ on $M$ define a Sasakian 3-structure if for even permutations $(i,j,k)$ of (1,2,3) the following properties are satisfied~:
\begin{eqnarray} \label{3-Sasakien}
\left\{ \begin{array}{rclcl}
	\varphi_k &=& \varphi_i \, \varphi_j - \omega_j \otimes R_i &=& - \varphi_j \, \varphi_i + \omega_i \otimes R_j \\
	R_k &=& \varphi_i (R_j) &=& - \varphi_j (R_i) \\
	\omega_k &=& \omega_i \circ \varphi_j &=& - \omega_j \circ \varphi_i,
\end{array} \right.
\end{eqnarray}
where $R_i$ is the Reeb vector field associated to $\omega_i$.

\begin{prop} --- \label{3-structures}
If a manifold $M$ of dimension $4n-1$ admits a Sasakian 3-structures given by three forms $\omega_1,\,\omega_2$ and $\omega_3$, then these forms generate a contact sphere which is both round and taut.
\end{prop}

{\bf Proof~:} Let $R_1,\,R_2$ and $R_3$ be the Reeb vector fields associated to $\omega_1,\,\omega_2$ and $\omega_3$. We define the tensor fields $\varphi_i$ by $g(X,\varphi_i(Y))=\frac{1}{2} \, d\omega_i(X,Y).$ Then for \mbox{$i=1,2,3$}, $(\varphi_i,\,R_i,\,\omega_i)$ defines a Sasakian structure and $\omega_1,\,\omega_2$ and $\omega_3$ define a Sasakian 3-structure. Therefore we have~:
\begin{eqnarray*}
\left\{ \begin{array}{rcl} 
g(X,\varphi_i(Y))&=&\frac{1}{2} \, d\omega_i(X,Y) \\ [1mm]
g(\varphi_i(X),\,\varphi_i(Y)) &=& g(X,Y)-\omega_i(X)\,\omega_i(Y) \\ [1mm]
\varphi_i^2 &=& -I + \omega_i \otimes R_i \\ [1mm]
\omega_i(Y) &=& g(R_i,\,Y),\; \textrm{for any} \; Y \\ [1mm]
\left[\varphi_i,\varphi_i \right] + d\omega_i \otimes R_i &=& 0,
\end{array} \right.
\end{eqnarray*}
for $i=1,\,2,\,3$ and (\ref{3-Sasakien}) is satisfied for even permutations $(i,j,\,k)$ of (1,2,3).

As $\omega_i(\varphi_i(X))=g(R_i,\varphi_i(X))=d\omega_i(R_i,X)=0$, we know that the image of $\varphi_i$ is in $Ker(\omega_i)$, for any $i$. Furthermore, as for any $X,Y \in Ker(\omega_i)$ we have $\varphi_i^2(X)=-X$ and $g(\varphi_i(X),\varphi_i(Y))=g(X,Y)$, $\varphi$ defines an isometry of $Ker(\omega_i)$ with $g(X,\varphi_i(X))=0$. For even permutations $(i,j,\,k)$ of (1,2,3) we also have $\omega_j(\varphi_i(X))=-\omega_k(X)$, so $\varphi$ preserves the intersection $\Sigma$ of the kernels of $\omega_1,\,\omega_2$ and $\omega_3$.

The Reeb vector fields $R_1,\,R_2$ and $R_3$ are everywhere linearly independant. Indeed, for even permutations $(i,j,\,k)$ of (1,2,3), we have~:
\begin{eqnarray}\label{orth}
\begin{array}{rclclclcl}
\omega_i(R_j) &=& \omega_j(\varphi_k (R_j)) &=& g(R_j,\varphi_k (R_j)) &=& \frac{1}{2} \, d\omega_k(R_j,R_j) &=& 0,
\end{array}
\end{eqnarray}
so $R_1,R_2$ and $R_3$ are orthogonal, as $\omega_i(R_j)=g(R_i,R_j)$.
None of these vectors is in $\Sigma$ and the dimension of $\Sigma$ is $4n-4$. Thus, for a given point $p \in M$, there is a basis of $T_p M$ which can be written as $(R_1,\,R_2,\,R_3,\,X_4,\dots,X_{4n-1})$, where the $X_i$ are elements of $\Sigma$. It is possible to choose them in such a way that for $k=4m$ we have $X_{k+1}=\varphi_1 (X_k), \quad X_{k+2}=\varphi_2 (X_k)$ and $X_{k+3}=\varphi_3 (X_k)$.

We set $\omega=\sum_{i=1}^3 \lambda_i \, \omega_i,$ with $\sum_{i=1}^3 \lambda_i^2=1$. Then we have 
$$\omega \wedge (d\omega)^{2n-1}(R_1,\,R_2,\,R_3,\,X_4,\dots,X_{4n-1}) = 2^{2n-2}\, (2n-2)!. $$
This is due to the equations
\begin{eqnarray*}
\omega_i(Y_j) &=& \delta_{ij}, \quad \textrm{for}  \; R_i = Y_i \; \textrm{and}  \; X_j = Y_j \\
d\omega_i(R_j,R_k) &=& 1, \;\; \quad \textrm{if} \; (i,j,\,k) \; \textrm{is an even permutation of} \; (1,2,3) \\
d\omega_i(R_j,X_k) &=& 0 \\
d\omega_i(X_j,X_k) &=& 0, \;\; \quad \textrm{if} \; X_k \neq \varphi_i\,X_j \\
d\omega_i(X_j,\varphi_i\,X_j) &=& -2.
\end{eqnarray*}
Thus, $\omega_1,\,\omega_2$ and $\omega_3$ generate a taut contact sphere.

Furthermore, we have (\ref{orth}) and for any vector field $X$
\begin{eqnarray*}
\begin{array}{crclcl}
R_i \lrcorner \, d\omega_j(X) &=& 2 \omega_k(X) &=& -R_j \lrcorner \, d\omega_i(X).
\end{array}
\end{eqnarray*}
So by Lemma \ref{CNS rond}, the contact sphere generated by $\omega_1,\,\omega_2$ and $\omega_3$ is round.
\finthm

This proposition yields new classes of examples of round and taut contact spheres in dimension higher than 3. According to Boyer, Galicki and Mann (see \cite{BGM}), the homogeneous manifolds which carry Sasakian 3-structures are the following~:
\begin{eqnarray*}
\frac{Sp(n+1)}{Sp(n)} \cong \sph^{4n+3}, \quad \frac{Sp(n+1)}{Sp(n) \times \mathbb{Z}_2} \cong \R P^{4n+3}, \\
\frac{SU(m)}{S(U(m-2) \times U(1))}, \quad \frac{SO(k)}{SO(k-4) \times Sp(1)}, \\
\frac{G2}{Sp(1)}, \quad \frac{F_4}{Sp(3)}, \quad \frac{E_6}{SU(6)}, \quad \frac{E_7}{Spin(12)}, \quad \frac{E_8}{E_7},
\end{eqnarray*}
for $m \geq 3$ and $k \geq 7$.
The named authors also prove the existence of infinitely many compact, simply connected and strongly inhomogeneous manifolds of dimension 7 which carry Sasakian 3-structures and which are homotopically distinct.
\vspace{4mm}

{\small
Acknowledgement: This paper is part of my PhD-Thesis (Mulhouse, december 2004). I am very grateful to my supervisors Robert Lutz and Norbert A'Campo for their constant support and many discussions about contact spheres and much more.}


{\small \begin{thebibliography}{33}
\bibitem{Ad} J.{\sc Adams}, Vector fields on spheres, Bull. Amer. Math. Soc., 68 (1962), 39-41.
\bibitem{Bl} D. {\sc Blair}, Contact Manifolds in Riemannian Geometry, Lecture Notes in Mathematics 509, Springer (1976).
\bibitem{BGM} C. {\sc Boyer}, K. {\sc Galicki}, B. {\sc Mann}, The geometry and topology of 3-Sasakian manifolds, J. reine u. angew. Math. 455 (1994), 183-220.
\bibitem{GG1} H.{\sc Geiges} and J.{\sc Gonzalo}, Contact geometry and complex surfaces, Invent. Math., 121 (1995), 147-209.
\bibitem{GG2} H.{\sc Geiges} and J.{\sc Gonzalo}, Contact Circles on 3-manifolds, J. Diff. Geometry, 46 (1997), 236-286.
\bibitem{Gr} J. W. {\sc Gray}, Some global properties of contact structures, Ann. of Math., 69 (1959), 421-450.
\bibitem{Eck} B.{\sc Eckmann}, Gruppentheoretischer Beweis des Satzes von Hurwitz-Radon über die Komposition quadratischer Formen, Comm. Math. Helv., 15 (1943), 358-366.
\bibitem{L1} R.{\sc Lutz}, Structures de contact sur les fibr\'es principaux en cercles de dimension trois, Ann. Inst. Fourier, Grenoble, 27, 3 (1977), 1-15.
\bibitem{L2} R.{\sc Lutz}, Sur la g\'eom\'etrie des structures de contact invariantes, Ann. Inst. Fourier, Grenoble, 29, 1 (1979), 283-306.
\bibitem{Ma} J.{\sc Martinet}, Sur les singularités des formes différentielles, Ann. Inst. Fourier, 20,1 (1970), 95-178.
\end{thebibliography}}
\end{document}